\documentclass{amsart}

\usepackage{amsmath,amsfonts,epsfig}

\newtheorem{theorem}{Theorem}[section]
\newtheorem{lemma}[theorem]{Lemma}
\newtheorem{proposition}[theorem]{Proposition}

\theoremstyle{definition}

\theoremstyle{remark}

\numberwithin{equation}{section}

\def\sv{\mathcal V}
\def\ghat{\hat{\gamma}}
\def\Ghat{\hat{\Gamma}}
\def\bbeta{\bar{\beta}}
\def\bgamma{\bar{\gamma}}

\begin{document}

\title{Relative hyperbolicity and right-angled Coxeter groups}

\author[Patrick Bahls]{Patrick Bahls}

\address{Department of Mathematics \\ University of Illinois at 
Urbana-Champaign \\ 1409 W. Green Street \\ Urbana, IL 61801 \\ USA}

\email{pbahls@math.uiuc.edu}



\begin{abstract}
We show that right-angled Coxeter groups are relatively hyperbolic in the 
sense defined by Farb, relative to a natural collection of rank-$2$ 
parabolic subgroups.
\end{abstract}

\thanks{The author was supported by an NSF VIGRE postdoctoral grant.}

\maketitle

\section{Introduction} \label{sectionintro}

The theory of (word) hyperbolic groups, initiated in the 1980s by M. 
Gromov in \cite{Gro1}, is a rich one.  The desirable geometric structure 
of such groups often provides a great deal of insight into the 
combinatorial structure of the group itself.  Gromov suggested a 
generalization of the notion of hyperbolicity, wherein a group would be 
said to be {\it relatively} hyperbolic with respect to a given collection 
of subgroups (known as parabolic subgroups).

In \cite{Farb}, B. Farb presented the first rigorous definition of 
relative hyperbolicity.  Though his definition is the one adopted in this 
paper, as well as in a number of others, it is not the only definition.  
B. Bowditch (\cite{Bow}) and A. Yaman (\cite{Yaman}) have presented 
alternative definitions of relative hyperbolicity.  A. Szczepa\'nski (in 
\cite{Sz1}) studied the relationship of these definitions, showing that 
the latter two are equivalent, and imply the first.  In \cite{Bow}, 
Bowditch shows that if a group is relatively hyperbolic in the sense of 
Farb and furthermore satisfies Farb's ``bounded coset penetration'' (BCP) 
property, then it is relatively hyperbolic in the sense of Bowditch and 
Gromov-Yaman.  If $G$ is relatively hyperbolic in the sense of Farb, it is 
often said to be relatively hyperbolic in the {\it weak} sense.

The objects of study in this paper are right-angled Coxeter groups.  
Recall that a {\it Coxeter system} $(W,S)$ is a pair in which $W$ is a 
group (called a {\it Coxeter group}) and $S=\{s_i\}_{i \in I}$ for which 
there is a presentation

$$\langle S \ | \ R \rangle$$

\noindent where

$$R = \{ (s_is_j)^{m_{ij}} \ | \ m_{ij} \in \{1,2,...,\infty\}, 
m_{ij}=m_{ji}, \ {\rm and} \ m_{ij}=1 \Leftrightarrow i=j \}.$$

\noindent (In case $m_{ij}=\infty$, the element $s_is_j$ has infinite 
order.)

Recall that for each subset $T \subseteq S$, the subgroup $W_T$ of $W$ 
generated by the elements of $T$ is also a Coxeter group, with system 
$(W_T,T)$ (see \cite{Bo}).  Such a subgroup is called a ({\it standard}) 
{\it parabolic subgroup}.

If $m_{ij} \in \{1,2,\infty\}$ for every $i,j \in I$, we call $(W,S)$ (and 
$W$) {\it right-angled}.  In his dissertation (\cite{Rad}), D. Radcliffe 
showed that the Coxeter presentation corresponding to a given right-angled 
Coxeter group is essentially unique, in that for any two fundamental 
generating sets $S$ and $S'$, there exists an automorphism $\alpha \in 
{\rm Aut}(W)$ such that $\alpha(S)=S'$.

The information contained in a Coxeter system can be captured graphically 
by means of a Coxeter diagram.  The {\it Coxeter diagram} $\sv=\sv(W,S)$ 
corresponding to the system $(W,S)$ is an edge-labeled graph with vertex 
set in one-to-one correspondence with the set $S$, and for which there is 
an edge $[s_is_j]$ labeled $m_{ij}$ between any two vertices $s_i \neq 
s_j$ whenever $m_{ij} < \infty$.  (In this paper, we omit the comma in the 
notation for an edge of a graph.)  It is easy to see that Radcliffe's 
rigidity result is equivalent to saying that to a given right-angled 
Coxeter group, there is a {\it unique} diagram, every edge of which has 
label $2$.  For this reason, we often suppress mention of the generating 
set $S$.

The following theorem follows from the work of G. Moussong 
(\cite{Moussong}):

\begin{theorem} {\bf [Moussong]} \label{theoremmoussong}
Let $W$ be a right-angled Coxeter group, with corresponding diagram $\sv$.  
Then $W$ is word hyperbolic if and only if $\sv$ contains no achordal 
simple circuits $\{[ab],[bc],[cd],[da]\}$ of length $4$.
\end{theorem}

That is, $W$ is hyperbolic provided we see no squares without diagonals in 
$\sv$.  (It is not hard to see that such squares imply the existence of 
flats in the Cayley graph of $W$, arising from the presence of the 
subgroup $D_{\infty} \times D_{\infty}$.)  A very natural way to attempt 
to overcome this obstacle to hyperbolicity would be to ``quotient'' by one 
of the diagonals of such a square, in effect collapsing one of the 
dimensions in the corresponding subgroup $H \cong D_{\infty} \times 
D_{\infty}$.  In fact, that is just what we will do.  The result is this 
paper's main theorem:

\begin{theorem} \label{theoremmain}
Let $W$ be a right-angled Coxeter group with diagram $\sv$.  For each 
achordal simple circuit $\{[a_ib_i],[b_ic_i],[c_id_i],[d_ia_i]\}$ in 
$\sv$, select a diagonal $\{a_i,c_i\}$.  Then $W$ is relatively hyperbolic 
(in the sense of Farb) relative to the collection of subgroups $\left\{ 
W_{\{a_i,c_i\}} \right\}$ for the above choice of diagonals.
\end{theorem}

\section{Relative hyperbolicity} \label{sectionrelhyp}

We now recall the definition of relative hyperbolicity due to Farb.

Suppose that $G$ is a finitely generated group with distinguished 
generating set $S$.  Let $\Gamma=\Gamma(G,S)$ be the Cayley graph for $G$ 
with respect to the set $S$.  We consider $\Gamma$ as a metric graph, 
where each edge has length $1$.  We will denote the vertex of $\Gamma$ 
corresponding to $g \in G$ by $v(g)$, or simply by $g$ when this notation 
will not be confusing.

Let $\{H_i\}_{i \in I}$ be a collection of subgroups of $G$ (for our 
purposes $|I| < \infty$).  We construct a new graph, 
$\Ghat=\Ghat(\{H_i\}_{i \in I})$ as follows.  For each left coset $gH_i$ 
($g \in G, i \in I$), we add a new vertex, denoted $v(gH_i)$, to $\Gamma$; 
for each element $g' \in G$ lying in the coset $gH_i$, we add an edge from 
the new vertex $v(gH_i)$ to the vertex $v(g') \in \Gamma$.  Define the 
length of each of the new edges to be $1/2$.  The graph $\Ghat$ so 
obtained is called the {\it coned-off Cayley graph of} $G$ {\it with 
respect to the collection} $\{H_i\}_{i \in I}$.

The effect of adding the new vertices is to create ``shortcuts'' which 
pass through left cosets; in particular, the distance from $v(1)$ to 
$v(g)$ is $1$ for any $g \in H_i$, for some $i \in I$.

Clearly the geometric structure of $\Ghat$ may be very different from that 
of $\Gamma$.  If $\Ghat$, with the metric defined in the paragraphs above, 
is a hyperbolic metric space, we say that $G$ is {\it relatively 
hyperbolic, relative to the collection} $\{H_i\}_{i \in I}$.  For details 
on hyperbolic metric spaces and their desirable properties, we refer the 
reader to \cite{BH}.  Here we merely remind that hyperbolic metric spaces 
can be roughly characterized by possessing ``thin'' geodesic triangles. 
(The precise condition that we will verify is contained below in 
Proposition~\ref{propositionpapa}.)

A number of results regarding relatively hyperbolic groups have been 
proven.  A good source of examples is the paper \cite{Sz2} by A. 
Szczepa\'nski.  Regarding the relative hyperbolicity of Coxeter groups and 
related groups, little has been shown.  In \cite{KapSch}, I. Kapovich and 
P. Schupp show that large-type Artin groups are relatively hyperbolic 
relative to a certain collection of rank-$2$ subgroups.  As was done in 
the latter paper, here we will rely upon the following useful result due 
to P. Papasoglu (in \cite{Papa}):

\begin{proposition} \label{propositionpapa}
Let $\Gamma$ be a connected graph with simplicial metric $d$.  Then 
$\Gamma$ is hyperbolic if and only if there is a number $\delta > 0$ such 
that for any points $x,y \in \Gamma$ (not necessarily vertices), any two 
geodesic paths from $x$ to $y$ in $\Gamma$ are $\delta$-Hausdorff close.
\end{proposition}

\noindent (See the remark in \cite{KapSch} regarding applicability of the 
result as originally stated to our situation.)

A graph with the simplicial metric can be obtained obtained from $\Ghat$ 
by subdividing each of the original edges of $\Gamma$.  However, it is 
not difficult to see that if the condition in 
Proposition~\ref{propositionpapa} can be established for the original 
metric on $\Ghat$, then the condition holds also for the simplicial metric 
(perhaps with a different constant $\delta$, of course).  In fact, it is 
clear that we need only verify that the condition holds for two {\it 
vertices} $x=v(g_1)$ and $y=v(g_2)$ in $\Gamma$.

\section{Admissible Coxeter groups}

We first focus our attention on a certain class of right-angled Coxeter 
groups which are (at least intuitively) easier to deal with than 
right-angled groups in general.  

Let $W$ be right-angled, and let $\sv$ be its diagram.  A {\it simple 
circuit of length $n$} in the diagram $\sv$ is a closed path 
$C=\{[s_1s_2],...,[s_ns_1]\}$ in $\sv$ such that $s_i \neq s_j$ for all 
$1 \leq i < j \leq n$.  We identify each circuit with its cyclic 
permutations.  A circuit $C$ as above is called {\it achordal} if for any 
$s_i,s_j$, $|j-i|>1$, $s_i$ and $s_j$ are not adjacent in $\sv$.

Any achordal simple circuit $\{ [ab],[bc],[cd],[da] \}$ in $\sv$ of length 
$4$ is called a {\it square}.  A {\it diagonal} of the square above is 
either one of the pairs $\{a,c\}$ or $\{b,d\}$. We call $\sv$ (and also 
$W$) {\it admissible} if it is possible to choose a collection of 
diagonals, one for each square in $\sv$, which are all disjoint from one 
another.  We will also call such a choice of diagonals admissible.  As 
examples, consider the Coxeter diagrams in Figure~1; $\sv_1$ represents 
an admissible group, and $\sv_2$, a group which is not admissible.  
(All edges are assumed to have label $2$.)

\begin{figure}[hbt]
	\begin{center}
		\input{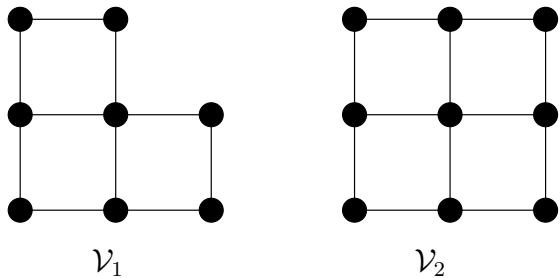}
	\end{center}
	\caption{Admissible and inadmissible groups}
\end{figure}

Now suppose that $\sv$ is admissible, and select an admissible set of 
diagonals, $\left\{ \{a_i,c_i\} \right\}$.  Let $H_i=W_{\{a_i,c_i\}} \cong 
D_{\infty}$.  We now prove the main theorem in the restricted setting of 
admissible Coxeter groups:

\begin{theorem} \label{theoremadm}
Let $W$ be an admissible Coxeter group.  Then $W$ is relatively hyperbolic 
relative to any collection of parabolic subgroups $\left\{ H_i \right\}$ 
chosen as above (that is, corresponding to an admissible choice of 
diagonals).
\end{theorem}

In order to prove this, we must verify that the condition in 
Proposition~\ref{propositionpapa} holds for admissible Coxeter groups.  We 
do this by explaining the structure of geodesic paths in $\Ghat=\Ghat 
\left( \{H_i\} \right)$ between two vertices $v(w_1)$ and $v(w_2)$ in 
$\Gamma=\Gamma(W,S)$.  Having described the structure of such paths, we 
will be able to show that any two such paths with the same endpoints are 
$\delta$-close, where $\delta$ is independent of the choice of endpoints.

Having proven Theorem~\ref{theoremadm}, we will (in the final section) 
indicate the minor changes that need be made in the arguments of 
Section~\ref{sectiongeods} in order to prove the main theorem in general.

\section{Geodesics in $\Ghat$} \label{sectiongeods}

Throughout this section, both the group $W$ and the choice of diagonals 
made below are assumed to be admissible.  We describe the collection of 
geodesics in $\Ghat$ from $1$ to a given vertex $w$ in $\Gamma$, in terms 
of one such geodesic.  Our description will enable an easy proof of 
Proposition~\ref{propositionpapa} in the presence of admissibility, and 
therefore of Theorem~\ref{theoremadm}.

First of all, let us note the following fact, which can be proven very 
easily using van Kampen diagrams.

\begin{lemma} \label{lemmasameletters}
Let $(W,S)$ be a right-angled Coxeter system, and let $w \in W$.  Any 
two geodesic words (in the generators $S$) representing $w$ contain the 
same number of occurrences of each letter $s \in S$.
\end{lemma}

Lemma~\ref{lemmasameletters} actually follows from a more general result 
(which can also be proven using van Kampen diagrams):

\begin{proposition} \label{propositionsameconjs}
Let $(W,S)$ be an arbitrary Coxeter system, and let $w \in W$.  Any two 
geodesic words (in the generators $S$) representing $w$ contain the same 
number of letter from each conjugacy class of generators from $S$.
\end{proposition}

Indeed, Lemma~\ref{lemmasameletters} follows from 
Proposition~\ref{propositionsameconjs} because in an even system (one in 
which every relator is of the form $(st)^m$, $m$ even), no two 
distinct generators are conjugate.  Lemma~\ref{lemmasameletters} will be 
used frequently without mention.

Now to the geodesics.  We will show that there are essentially two 
operations which can be applied in order to obtain new geodesics from a 
fixed geodesic.  Having described these operations, we will show that the 
geodesics obtained from a single geodesic by means of these operations are 
uniformly close to the original geodesic.

Let $\Gamma=\Gamma(W,S)$ and $\Ghat=\Ghat \left( \{H_i\} \right)$ as in 
the preceding section.  Let $v(w_1)$ and $v(w_2)$ be vertices in $\Gamma$.  
By translation, we can assume that $w_1=1$, and we rewrite $w_2$ as $w$.  
Through abuse of notation, we will often write the vertex $v(w)$ as $w$.

Suppose that $\ghat$ is a $\Ghat$-geodesic from $1$ to $w$.  The path 
$\ghat$ can be broken up into blocks which correspond to subpaths in 
$\Gamma$ and subpaths in $\Ghat \setminus \Gamma$.  That is, 
$\ghat=\alpha_1 \beta_1 \cdots \alpha_k \beta_k$, where $\alpha_i$ is a 
path

$$\left\{ [w_{i-1},w_{i-1}s_1],[w_{i-1}s_1,w_{i-1}s_1s_2], ... 
[w_{i-1}s_1s_2 \cdots s_l,w_i] \right\}$$

\noindent where every vertex and every edge lies in $\Gamma$, and 
$\beta_i$ is a path of combinatorial length $2$:

$$\left\{ [w_{i-1},v(w_{i-1}H_i)],[v(w_{i-1}H_i),w_i=w_{i-1}h_i] 
\right\},$$

\noindent for some $h_i \in H_i$.  We call the paths $\alpha_i$ the 
$\Gamma$-{\it blocks} of $\ghat$ and the paths $\beta_i$ the $\Ghat$-{\it 
blocks} of $\ghat$.

One fact is obvious (by the geodesity of $\ghat$):

\begin{lemma} \label{lemmaalphageod}
Each $\Gamma$-block $\alpha_i$ of $\ghat$ is a geodesic in $\Gamma$.
\end{lemma}

When it will not cause confusion, we do not distinguish between the path 
$\alpha_i$ and the word which labels this path.

We now construct a path, $\gamma$, in $\Gamma$, from $\ghat$.  For each 
$\Ghat$-block $\beta_i$ in $\ghat$ as above, replace $\beta_i$ with the 
path $\bbeta_i$ in $\Gamma$ from $w_{i-1}$ to $w_i=w_{i-1}h_i$ defined by 
the (unique) geodesic word in $W_{\{a_i,c_i\}}$ representing $h_i$.  If 
$\bbeta_i$ is a single letter, we call the corresponding $\Ghat$-block 
{\it trivial}; otherwise, we call the block {\it non-trivial}.  We will 
often abuse terminology and call $\bbeta_i$ (or even the label of this 
path) a $\Ghat$-block as well.

Performing all such replacements $\beta_i \mapsto \bbeta_i$ yields the 
path $\bgamma$.  In general, $\bgamma$ may not be a $\Gamma$-geodesic, but 
it is nearly so.  Namely, we have the following

\begin{lemma} \label{lemmagammageod}
Let $\bgamma$ be derived from $\ghat$ as above.  Let $h_i$ be the geodesic 
word in $H_i$ corresponding to the $\Ghat$-block $\beta_i$.  Then a 
$\Gamma$-geodesic $\gamma$ from $1$ to $w$ can be obtained from $\bgamma$ 
by canceling at most $2$ letters from each word $h_i$ arising from a 
non-trivial $\Ghat$-block.
\end{lemma}

\begin{proof}
We appeal to an algorithm of Tits in \cite{Tits}.  The results of this 
paper show that if a word $w$ in the letters $S$ is not geodesic, then a 
shorter representative for the same group element can be obtained from $w$
by performing successive commutations to bring two occurrences of the same 
letter next to each other, which can then be canceled.

It is clear that there can be no such cancellation of the letter $a$, 
where either $a$ occurs in both $\alpha_{i_1}$ and $\alpha_{i_2}$ or $a$ 
occurs in both $\alpha_{i_1}$ and $h_{i_2}$.  (Otherwise we would be
able to replace the appropriate words with shorter words and readjust the 
path $\bgamma$ to obtain a shorter (relative to both both $\Gamma$ and 
$\Ghat$) path from $1$ to $w$, contrary to the geodesity of $\ghat$.)  The 
same holds for trivial $\Ghat$-blocks $\beta_i$.

Thus we need only consider the case where the letter $a$ to be canceled 
occurs in two words $h_{i_1}$ and $h_{i_2}$ ($i_1<i_2$) coming from 
non-trivial $\Ghat$-blocks $\beta_{i_1}$ and $\beta_{i_2}$.  (By 
admissibility, the other letter, $c$, in $h_{i_1}$ is the same as the 
other letter in $h_{i_2}$.)  Suppose cancellation of $a$ does occur.  In 
this case, since $ac \neq ca$, $h_{i_1}$ ends with $a$ and $h_{i_2}$ 
begins with $a$, and $a$ commutes with every letter occurring in any word 
$\alpha_j$ or $h_j$ between $h_{i_1}$ and $h_{i_2}$.  The letter $c$ 
cannot have such commutativity properties, as otherwise we would be able 
to commute the entire word $h_{i_1}$ past the intervening words and 
collapse the two $\Ghat$-blocks $\beta_{i_1}$ and $\beta_{i_2}$.  We would 
thus obtain a $\Ghat$-path from $1$ to $w$ with subpaths $\alpha_i$ 
identical to those of $\bgamma$, but with one fewer $\Ghat$-block, 
contradicting $\ghat$'s geodesity.  Therefore we may cancel only one 
letter from the end of $h_{i_1}$ and one from the beginning of $h_{i_2}$.

It may be possible (if $h_i$ has length at least $3$) to cancel both the 
first and the last letters of $h_i$ in this fashion.
\end{proof}

\noindent {\bf Remark.} \ One can use van Kampen diagrams to prove 
the above result, instead of appealing to \cite{Tits}.  For applications 
of van Kampen diagrams to Coxeter groups and related groups, see 
\cite{Ba}, \cite{BaMi}, and \cite{KapSch}, for example.  For details 
regarding van Kampen diagrams, the reader may consult \cite{LynSch} or 
\cite{Ols}.

\vskip 1mm

The proof of Lemma~\ref{lemmagammageod} suggests a certain rigidity 
property that $\Ghat$-geodesics must satisfy.  We say that two 
$\Ghat$-blocks $\beta_1$ and $\beta_2$ have the same {\it type} if they 
are both words in letters of the same diagonal, $\{a,c\}$.  Arguments 
similar to the proof above can be used to verify

\begin{lemma} \label{lemmatypeofbetas}
Let $\ghat_1$ and $\ghat_2$ be two $\Ghat$-geodesics from $1$ to $w$.  
Then the non-trivial $\Ghat$-blocks of $\ghat_1$ are the same type, up to 
multiplicity, as those of $\ghat_2$.  (In fact, corresponding 
$\Ghat$-blocks $\beta$ and $\beta'$ give rise to words $h$ and $h'$ which 
differ by at most $2$ letters, as in Lemma~\ref{lemmagammageod}.)
\end{lemma}

As a consequence of these lemmata, there is only one fundamental operation 
by means of which we may obtain the $\Ghat$-blocks of one 
$\Ghat$-geodesic, $\ghat_1$, from the $\Ghat$-blocks of another, 
$\ghat_2$, with the same endpoints.  Namely, suppose $\beta_{i_1}$ and 
$\beta_{i_2}$ are of the same type (corresponding to $H=W_{\{a,c\}}$).  
Suppose that the letter $a$ commutes with every letter in any block lying 
between $\beta_{i_1}$ and $\beta_{i_2}$.  Our operation then consists of 
replacing $h_{i_1}$ with $h_{i_1}a$ and $h_{i_2}$ with $ah_{i_2}$.  We 
call an application of this operation an {\it insertion-deletion} (ID).  
ID operations may be performed a number of times on a given geodesic.

Now note that it is possible that some $\Ghat$-blocks may commute with one 
another and so the relative order of the blocks in the geodesics $\ghat_1$ 
and $\ghat_2$ may not be the same.  Likewise, given two $\Ghat$-geodesics, 
the letters which appear in the $\Gamma$-blocks of each geodesic may 
appear in a different order in each path as well.  Furthermore, some such 
letters may commute with certain of the $\Ghat$-blocks, leading to greater 
variation.  However, this variation is controllable: we now define a 
procedure by means of which every $\Ghat$-geodesic from $1$ to $w$ can be 
obtained, once we have in hand one such geodesic, $\ghat$.

Let

$$\bgamma = \alpha_1 \bbeta_1 \cdots \alpha_k \bbeta_k$$

\noindent be obtained from $\ghat$ by means of the replacements $\beta_i 
\mapsto \bbeta_i$.  As before, let $h_i$ be the unique geodesic word 
labeling $\bbeta_i$.  Let $w = \alpha_1 h_1 \cdots \alpha_k h_k$; that is, 
$w$ is the word labeling $\bgamma$.

We now describe all words $\bgamma_1$ (with label $w_1 =_W w$) which can 
be obtained from a $\Ghat$-geodesic $\ghat_1$ which in turn has the same 
endpoints as $\ghat$.  Such a word $\bgamma_1$ will be called {\it nice}; 
its label will also be called nice.  In particular, each subword $h_i$ of 
$w$ appears undisturbed in the label of any nice word.  (As an example, if 
$\sv$ consists of the square $\{ [ab],[bc],[cd],[da] \}$, and if $\{a,c\}$ 
is the selected diagonal, then $acab$ and $acba$ are both geodesics 
representing the same element as $baca$, whereas only the first is nice.)  
We assume (often without mention) throughout the remainder of the section 
that every word $\bgamma$ is nice.

We define a {\it syllable} of $w$ to be either a letter in some 
$\Gamma$-block $\alpha_i$ or a single $\Ghat$-block $h_i$ in its entirety.
We write $\|w\|$ for the {\it syllable length} of $w$ (that is, the number 
of syllables occurring in $w$).  Clearly every nice word $w_1$ satisfies 
$\|w\|=\|w_1\|$.

Certain of the syllables of the word $w$ are forced to occur before 
certain others, in {\it any} nice word $w_1$.  That is, for instance, if 
$st \neq ts$ and the first occurrence of $s$ precedes the first occurrence 
of $t$ in $w$, then in any nice word $w_1$, the first occurrence of $s$ 
must precede the first occurrence of $t$.  It is possible (although not 
necessary for our purposes) to assemble a complete description of these 
precedence relationships that arise in the word $w$.

\vskip 1mm

\noindent {\bf Example.} \ Consider the group $W$ with diagram shown in 
Figure~2.

\begin{figure}[hbt]
	\begin{center}
		\input{fig2.tex}
	\end{center}
	\caption{}
\end{figure}

We select the diagonals $\{a,c\}$ and $\{b,f\}$.  Consider the word 
$acabedbfbc$ (which has $6$ syllables).  In any nice word $w_1$ such that 
$w =_W w_1$, $aca$ must precede each of $e$, $d$, $c$, and $bfb$, but may 
follow the first $b$.  Also, $d$ must follow $aca$, $b$, and $e$, and must 
precede $bfb$, but may either precede or follow $c$.

\vskip 1mm

From similar considerations, it is possible to determine, for a given 
syllable $x$ of $w$, which syllables of $w$ may appear {\it before} $x$ in 
some nice word $w_1$ such that $w =_W w_1$.  Of course, any syllable 
preceding $x$ in $w$ may precede $x$ in some $w_1$.  Furthermore, any 
syllable $x'$ which occurs after $x$ in $w$ and which commutes with both 
every syllable lying between $x$ and $x'$, {\it and} with $x$ itself can 
{\it also} precede $x$ in some nice word $w_1$.

\vskip 1mm

\noindent {\bf Example.} \ Returning to the previous example, $d$ can be 
preceded by each of $aca$, $b$, $e$, and $c$.

\vskip 1mm

\noindent {\bf Remark.} \ As the reader may have noticed, for the purposes 
of establishing rules of precedence, as above, we may consider distinct 
occurrences of the same letter $s$ as different letters (or as indexed 
copies of the same letter).  This fact is used implicitly throughout the 
remainder of this section.

\vskip 1mm

Consider any syllable $x$ occurring in $w$.  We say that $x$ is $k$-{\it 
forced} if in any nice word $w_1 =_W w$, $x$ must appear in the first $k$ 
syllables of $w_1$.  For instance, every syllable in $w$ is 
$\|w\|$-forced, and in the above example, the syllable $d$ is $5$-forced.

As we will soon see, any prefix of a nice word $w_1$ of syllable length 
$k$ will contain a certain number of $k$-forced syllables, as well 
(perhaps) as a number of other syllables are not $k$-forced.  Syllables 
that occur in a prefix of $w_1$ of syllable length $k$ which are not 
$k$-forced will be called $k$-{\it free}.

\begin{lemma} \label{lemmaforcedfree}
Let $w_1$ be a nice word such that $w =_W w_1$, and let $p(k)$ be the 
prefix of $w_1$ of syllable length $k$.  Then $p(k) =_W p_1(k)p_2(k)$, 
where every syllable in $p_1(k)$ is $k$-forced, and every syllable in 
$p_2(k)$ is $k$-free.
\end{lemma}

\begin{proof}
Let $x$ be any $k$-free syllable in $p(k)$, and let $y$ be a $k$-forced 
syllable occurring after $x$ in $p(k)$.  It must be possible to commute 
$x$ past $y$ (and indeed past every intervening syllable), as otherwise 
$x$ too would be $k$-forced.  Therefore every $k$-free syllable can be 
commuted past every $k$-forced syllable, and $p(k)$ can be rewritten as 
claimed.
\end{proof}

Being $k$-forced does not depend on the choice of nice word, so the word 
$p_1(k)$ contains the same syllables, regardless of our choice of nice 
word $w_1$.  Moreover, although these syllables are bound by rules of 
precedence (as above), it is clear that the word $p_1(k)$ obtained from a 
nice word $w_1$ can be modified by commutations to obtain the word 
$p_1(k)$ corresponding to any other nice word $w_2$.  The following fact 
results:

\begin{lemma} \label{lemmakforcedword}
We may choose the same word $p_1(k)$ (as in Lemma~\ref{lemmaforcedfree}) 
for every nice word $w_1$ such that $w =_W w_1$.
\end{lemma}

This fact is promising, as it indicates how we may prove that any two 
$\Ghat$-geodesics are uniformly close to one another.  Indeed, suppose
that for every $k$, $k-\|p_1(k)\| \leq M$ for some $M$, and let 
$\gamma(p_1)$ be the $\Ghat$-geodesic from $1$ to $v=v(p_1)$ in $\Ghat$ 
which is determined by $p_1$ in the obvious manner.  Let $w_1$ and $w_2$ 
be nice words such that $w =_W w_i$, $i=1,2$, and such that no ID 
operations need be performed in transforming $w_1$ to $w_2$.  Let $v_i$ 
be the point on the path in $\Ghat$ with label $w_i$ of $\Ghat$-distance 
$k$ from $1$.  Then

$$d_{\Ghat}(v_1,v_2) \leq d_{\Ghat}(v_1,v) + d_{\Ghat}(v,v_2) \leq 2M.$$

If we can uniformly bound $M$ (that is, show that it is independent of the 
choice of $w$), we will be nearly done.  (Because in this case, 
corresponding points on geodesic paths obtained from one another without 
the use of ID operations will be uniformly close.)

For a fixed system $(W,S)$ and admissible choice of diagonals, let $M$ be 
one less than the maximal number of syllables which mutually commute with 
one another in any possible word $w$.  (For instance, in the above 
examples, $M=1$.)  The value of $M$ can be determined easily from the 
diagram $\sv$ and from the choice of diagonals.

\begin{proposition} \label{propositionM}
Let $w$ label a nice $\Gamma$-geodesic $\gamma$.  For each $k$, $1 \leq k 
\leq \|w\|$, define $p_1(k)$ as above.  Then for every $k$, $1 \leq k \leq 
\|w\|$, $k-\|p_1(k)\| \leq M$.
\end{proposition}

\begin{proof}
For any word $w$ and any $k$, $1 \leq k \leq \|w\|$, denote by $n(k,w)$ 
the number of syllables of $w$ which are $k$-forced in $w$.  We need to 
show that $k-n(k,w) \leq M$, for every $k$, $1 \leq k \leq r$.

We prove this fact by induction on the syllable length of $w$.  The result 
is clearly true if $\|w\| = 1$.  Assume it to be proven in case $\|w\| 
\leq r-1$, and let $\|w\| = r$.

From $w$ we create a shorter word, $w'$, by removing a single syllable.  
Namely, let $x$ be any syllable of $w$ which is $k$-forced, for $k$ 
minimal among all syllables of $x$.  We claim that $x$ may be commuted to 
the front of $w$.  If this were not the case, some syllable, $y$, which 
precedes $x$, satisfies $xy \neq yx$.  In this case, strictly fewer 
syllables of $w$ can precede $y$ as can precede $x$, so that $y$ is 
$k$-forced for a smaller value of $k$ than is $x$, contrary to our choice 
of $x$.

Therefore $w =_W xw'$ for some $w'$, $\|w'\|=r-1$.  It is clear that $w'$ 
is a nice word (that is, that it labels a $\Gamma$-geodesic $\gamma'$ 
which arises from a $\Ghat$-geodesic $\ghat'$ by replacements $\beta_i 
\mapsto \bbeta_i$ as before).

Consider any syllable $y$, $y \neq x$ in $w$ (the same syllable appears in 
$w'$).  If $y$ is $k$-forced in $w'$, then it is clearly $(k+1)$-forced in 
$w$ ($x$ may precede it).  Moreover, $x$ is $(k+1)$-forced in $w$ as well.  
Thus $n(k+1,w) \geq n(k,w')+1$.

Fix $k \geq 2$ such that $n(k,w) \geq 2$.  By the inductive hypothesis, 

$$(k-1)-n(k-1,w') \leq M.$$

\noindent Thus

$$k-n(k,w) \leq k-(n(k-1,w')+1) = (k-1)-n(k-1,w') \leq M.$$

Therefore we need only show that $k-n(k,w) \leq M$ when $n(k,w)=1$; that 
is, when $x$ is the only syllable which is $k$-forced.  It suffices to 
show that for any choice of $w$, there is some syllable which is 
$(M+1)$-forced.

As above, let $x$ denote any syllable which is $k$-forced, for $k$ minimal 
among all $x$ in $w$.

We construct a set $S_1$ of syllables of $w$, all of which commute with 
one another, in the following fashion.  Initially, $S_1$ consists of the 
first syllable of $w$.  We add to this set the second syllable of $w$ if 
it commutes with the first.  Thereafter, we consider each syllable $y$ in 
turn, adding it to $S_1$ if and only if

\vskip 1mm

\noindent 1. \ $y$ commutes with each syllable already included in $S_1$, 
and

\vskip 1mm

\noindent 2. \ $y$ can be commuted to the front of the word $w$.

\vskip 1mm

\noindent It is easy to show that $x \in S_1$.

Let $m=|S_1|$.  By the definition of $M$, $m \leq M+1$.  We now claim that 
$x$ is $m$-forced (this will finish our proof).

Suppose that $S_1=\{x_1,...,x_{m-1},x\}$.  Without loss of generality, we 
can write $w$ as $x_1 \cdots x_{m-1}x w''$ for some word $w''$.

Assume that $x$ is not $m$-forced, and assume first that $x$ is a single 
letter.  Thus there is some syllable $y \not\in S_1$ such that $yx=xy$, 
but $yx_i \neq x_iy$, for some $i$, $1 \leq i \leq m-1$.  We may assume 
that $y$ is the first such syllable to occur after $x$ in $w$, and by 
commuting, we may assume that $i=m-1$.  Assume for now that both $x_{m-1}$ 
and $y$ are also single letters.  By the minimality in our choice of $x$, 
there must be a syllable $y'$ following $x$ in $w$ such that $x_{m-1}$ and 
$y$ both commute with $y'$, but $x$ does not.  (In fact, from our choice 
of $y$, this must be true of the first syllable following $x$ which does 
not commute with $x$.)  Assume for now that $y'$ too is a single letter.

We have found a square in the diagram $\sv$, consisting of the letters 
$\{x,x_{m-1},y',y\}$.  Therefore either $\{x,y'\}$ or $\{x_{m-1},y\}$ must 
have been selected as a diagonal.  However, either choice violates the 
assumption that each of the letters $x$, $x_{m-1}$, $y$, and $y'$ is a 
syllable in its own right.  (This is so because in either case we realize 
a non-trivial $\Ghat$-block consisting of two of these letters.)

There are a number ($15$, to be exact) of other cases to consider, 
depending on which of the syllables $x$, $x_{m-1}$, $y$, and $y'$ are 
single letters, and which stem from non-trivial $\Ghat$-blocks.  In each 
case, we may argue much as above, obtaining either a similar 
contradiction, or a contradiction to the very admissibility of the group 
$W$.  (The latter contradiction arises, for instance, when $x$ and $y$ are 
single letters and $x_{m-1}$ and $y'$ come from non-trivial 
$\Ghat$-blocks.)

Therefore $x$ is $m$-forced, and is thus $(M+1)$-forced.  This concludes 
the proof.
\end{proof}

What have we now shown?  Let $\ghat_1$ and $\ghat_2$ both be 
$\Ghat$-geodesics between the vertices $1$ and $w$ in $\Ghat$.  Then if 
$\ghat_2$ can be obtained from $\ghat_1$ without applying ID operations 
({\it i.e.}, by commutations of syllables only), $\ghat_1$ and $\ghat_2$ 
are $2M$-Hausdorff close, where $M$ is defined as above.

In order to prove Proposition~\ref{propositionpapa} (and therefore 
Theorem~\ref{theoremadm}), it suffices to show that if $\ghat_2$ and can 
be obtained from $\ghat_1$ by applying {\it only} ID operations, then the 
paths are $\delta$-close, for some $\delta > 0$.  (We can first apply 
ID operations to make sure that all $\Ghat$-blocks which occur in two 
given nice geodesics are identical up to commutation, and then perform 
commutations to obtain one geodesic from the other.)

Indeed, this is so.  The $\Ghat$-distance between corresponding points on 
two $\Ghat$-geodesics which differ only by application of ID operations is 
uniformly bounded by a constant which (like $M$) is related to the 
structure of $\sv$ and the choice of diagonals.

Each ID operation, like a balanced pair of parentheses, has an ``opening'' 
and a ``closing''.  In the notation used before to define ID operations, 
the ``opening'' occurs when the word $h_{i_1}$ is replaced by the word 
$h_{i_1}a$, and the closing when $h_{i_2}$ is replaced by $ah_{i_2}$.  The 
greatest distance between corresponding points on two geodesics differing 
only by ID operations is obtained when the number of ``unclosed'' pairs is 
maximal.  (This distance at most the number of unclosed pairs at that 
point).  A necessary (though not in general sufficient) condition for each 
of the ID operations inserting and deleting the letters $a_1,a_2,...,a_k$ 
to be simultaneously in an open state is that each $a_i$ commutes with  
$a_j$, $j>i$.  Thus the maximal number of unclosed pairs at any time is 
bounded above by the maximal number of generators which commute with a 
given generator of $S$.

Therefore, Proposition~\ref{propositionpapa} and Theorem~\ref{theoremadm} 
are proven.  Moreover, an explicit ``constant of relative hyperbolicity'' 
can be computed, merely by examining the diagram $\sv$ and the admissible 
choice of diagonals.

\section{The main theorem in general} \label{sectionmain}

We now indicate the changes that must be made in the above argument in 
order to prove Theorem~\ref{theoremmain} in general.  The chief difficulty 
stems from the fact that if $W$ (or the choice of diagonals) is not 
admissible, then a single letter $s \in S$ may appear in more than one 
diagonal.

However, the first four results of Section~\ref{sectiongeods} remain true 
in the general case.  In particular, we need only cancel at most $2$ 
letters from each non-trivial $\Ghat$-block in order to obtain a 
$\Gamma$-geodesic from a $\Ghat$-geodesic.  Moreover, ID operations are 
still valid, although now an ID operation can effect the exchange of a 
letter between two $\Ghat$-blocks of different type (which share a single 
letter).

Lemma~\ref{lemmatypeofbetas} is no longer valid.  As before, given two 
$\Ghat$-geodesics with the same endpoints in $\Ghat$, $\ghat_1$ and 
$\ghat_2$, we obtain $\Gamma$-geodesics $\gamma_1$ and $\gamma_2$ by first 
replacing each $\beta_i$ with $\bbeta_i$, and then by applying ID 
operations.  In order to obtain $\gamma_2$ from $\gamma_1$, we must modify 
the $\Ghat$-blocks that occur in $\gamma_1$ without changing the syllable 
length of its label.  This can be done step-by-step, where at each 
step we either leave fixed the number of $\Ghat$-blocks (merely shifting 
letters from one block to another), or decrease the number of 
$\Ghat$-blocks by one, adding $1$ to the total length of the 
$\Gamma$-blocks $\alpha_i$.

An operation which accomplishes this first goal is called a {\it shift} 
operation.  Let $\{a,b\}$, $\{b,c\}$, $\{c,d\}$, and $\{d,e\}$ all be 
selected diagonals.  Suppose that the word $w$ contains the subword 
$(ab)^mw_1cdw_2(ed)^n$, and that both $bw_1=w_1b$ and $dw_2=w_2d$ hold in 
$W$.  The most general form of shift involves replacing the subword above 
with the subword $(ab)^{m-1}aw_1bcw_2(de)^nd$ (both of which have the same 
number of syllables), or {\it vice versa}.  Of course, the words $w_i$ may 
be trivial.  More than one shift can be performed in an overlapping 
sequence, each syllable passing on a single letter to the next, in the 
obvious fashion.  In this case only the first and the last syllables 
involved in the shift can have length greater than $2$.

A second way in which the first goal can be accomplished is by applying an 
{\it exchange} operation.  Suppose that $\{a,b\}$, $\{a,c\}$, $\{b,d\}$, 
and $\{c,d\}$ are all selected diagonals, and that $bc=cb$.  Furthermore, 
assume that both $b$ and $c$ commute with the word $w_1$.  In its most 
general form, an exchange consists of replacing the subword $abw_1cd$ of 
$w$ with the subword $acw_1bd$.  As with shifts, there can be sequences of 
overlapping exchanges.

How do we effect a transformation in which a $\Ghat$-block disappears and 
the total length of the $\Gamma$-blocks $\alpha_i$ is increased by $1$?  
The disappearing $\Ghat$-block must have length $2$, and one of its 
letters must be ``absorbed'' into another $\Ghat$-block.  We've seen this 
behavior already, in ID operations which replace $(ab)^mw_1(ac)$ with 
$(ab)^maw_1c$, for example.  Another way this can be accomplished is by a 
variant sort of exchange, in which the word $(ab)w_1(cd)w_2(ed)^n$ is 
replaced by the word $(ac)bw_1w_2(de)^nd$ (or {\it vice versa}), where 
$\Ghat$-blocks are indicated by parentheses.

These new operations differ from ID operations in that the type of 
syllables present is (in general) modified.  However, each single new 
operation can be seen as a composition of two ID operations, where the 
word resulting from the application of the first ID operation is not a 
$\Ghat$-geodesic.  (In an exchange, there is a commutation of letters 
in between the two ID operations.)  For this reason, we may treat 
shifts and exchanges much as we would ID operations.  Indeed, keeping 
track of ``openings'' and ``closings'' of sequences of shifts and 
exchanges as we did with ID operations, it is not difficult to verify the 
following fact.

\begin{proposition} \label{propositionIDshiftexch}
Let $\ghat_1$ and $\ghat_2$ be $\Ghat$-geodesics which can be obtained 
from one another through application of ID operations, shifts, and 
exchanges only.  Then corresponding points on these two geodesics are 
uniformly $\Ghat$-close.
\end{proposition}

We are almost done now.  We first note the following analogue of 
Lemma~\ref{lemmatypeofbetas}:

\begin{proposition} \label{propositionthreeops}
Let $\ghat_1$ and $\ghat_2$ be $\Ghat$-geodesics with the same endpoints, 
$1$ and $w$, in $\Ghat$.  Then the $\Ghat$-blocks of $\ghat_2$ can be 
obtained from the $\Ghat$-blocks of $\ghat_1$ by application of ID 
operations, shifts, and exchanges.
\end{proposition}

\begin{proof}
Let $\bgamma_i$ be the $\Gamma$-path resulting from $\ghat_i$ be the 
replacements $\beta_i \mapsto \bbeta_i$ as before.  Both $\bgamma_1$ and 
$\bgamma_2$ are nice words representing the same group element.  By 
applying ID operations, we can transform each of these paths into a 
$\Gamma$-geodesic.  The label of one of these paths can be obtained from 
the label of the other by commutations whose net effect must not increase 
the syllable length of the word in question.  We have defined the 
operations above precisely in order to describe all ways in which this can 
be done.
\end{proof}

Now suppose we are given two $\Ghat$ geodesics, $\ghat_1$ and $\ghat_2$, 
with the same endpoints, $1$ and $w$, in $\Gamma$.  By applying ID 
operations, shifts, and exchanges, we may obtain from $\ghat_1$ a new 
$\Ghat$-geodesic, $\ghat'_1$, which has the same $\Ghat$-blocks as 
$\ghat_2$, although in a different order.  By 
Proposition~\ref{propositionIDshiftexch}, corresponding points on  
$\ghat_1$ and $\ghat'_1$ are uniformly close to one another in $\Ghat$.  
In order to prove Theorem~\ref{theoremmain}, we need only show that 
geodesics which differ only by the order of their syllables are uniformly 
close.

However, the proof of Proposition~\ref{propositionM} in the more general 
setting goes through almost exactly as before.  The only difference is 
that now we can no longer derive a contradiction to admissibility (such 
admissibility assumptions are not made!).  We only contradict the 
syllabification of the geodesic, as was done in the proof of 
Proposition~\ref{propositionM} above.

\end{document}